# Dynamic Model of Smoothing Problem in Water Power Systems

Jimsher Giorgobiani
Mziana Nachkebia
Weldon A. Lodwick

**ABSTRACT. In this paper the problem of optimal performance of a power system is considered. The problem is posed in various aspects within the frames of the theory of optimal control of stores. Mathematical models are presented by means of the recurrent equations of a dynamic programming. In the general case a new method is presented which avoids the "curse of dimensionality."**

The relationship of operational hydroelectric power stations and systems of water reservoirs is studied in this research with the objective of controlling water stores in order to satisfy demand for the electric power with the least expenditures. For any decision period, the share of hydroelectric energy with respect to total energy demand must be established. Hydroelectric energy is assumed to be free-of-charge in comparison with expensive thermal energy. Moreover, because of the specific character of river flow and demand of electric power the problem of providing energy at the least cost is dynamic and stochastic.

Many authors, among them [4], investigated the problem of optimal performance of a power system with large number of various power plants by means of the models of mathematical programming. For a single hydroelectric power station this problem under different assumptions was investigated by method of dynamic programming by J. Little [1], S. Karlin [2], J. Gessford [3]. They investigated a problem from the point of view of the theory of optimal control of stores. In this study, we apply the refined method of functional equations to investigate more general problems and to remove some assumptions imposed in the papers cited above.

The planning and analysis horizon of scheduling (year, as a rule) is divided into $n$ periods (fortnight, monthly, etc.); during one period a flow of water, production and demand for the electric power are considered uniform in time.

We introduce the following notations:

$x_i$ - inflow of water in water reservoir during the period $i$,

$u_i$ - outflow of water from water reservoir to turbines during the period $i$,

$r_i$ - demand for the electric power during the period $i$,

$c$ - price of the unit of electric power from thermal stations,

$K$ - maximal energy produced by thermal station during a period,

$p$ - penalty for deficiency per unit of the electric power ($p > c$),

$Q_i$ - store of water in water reservoir at the beginning of $i$-th period ($Q_i \leq \overline{Q}$).

Parameters $x_i, u_i, Q_i$ can be understood both as an amount of water, and an amount of energy as well, therefore it is possible to measure them in terms of volume (cubic meters) as well as in terms of the electric power (kilowatt-hours). Conversion of the amount of water to the corresponding amount of energy is performed using the well known formula, $N = 9.8\eta u H$, where $\eta$ - is the efficiency factor and $\overline{H}$ is the maximal flow where flow $H = H(Q) \leq \overline{H}$ is the function of the amount of water in water reservoir and varies in time.

For all existing in G water reservoirs in a given problem, power engineers, using the topographic maps, have established the relations $H = H(Q)$ by means of curves. Applying these graphs for our problem it is possible to design tables with appropriate step in advance and calculate the energy $N_i$ for period $i$ in kilowatt-hours by the formula

$$N_i = 9.8\eta \sum_{t=1}^{\tau} H\left(\min\left\{Q_i + t \cdot \frac{x_i - u_i}{\tau}, \overline{Q}\right\}\right) \cdot \frac{u_i}{\tau},$$

where $\tau$ is the number of days in a period.

Below all parameters are measured by equivalent amounts of the electric power. So we use $u_i$ instead of $N_i$. Moreover, we'll consider a base case and create for this case a recurrent set of equations of a dynamic programming. Thus, the expenditure during one period (the $i^{th}$ one) will be

$$c \cdot \min\{K, r_i - u_i\} + p \cdot \max\{0, r_i - u_i - K\}.$$

If we denote the minimum total expenditures for the future periods from $i$ up to and including $n$ (the time for which we are analyzing the optimal control of water store when there is $Q_i$ amount of water in a water reservoir at the beginning of $i^{th}$ period), by $f_i(Q_i)$, then the following recurrent relation for $i = 1, 2, \ldots, n-1$ is obtained:

$$f_i(Q_i) = \min_{u_i}\left\{c \cdot \min\{K, r_i - u_i\} + p \cdot \max\{0, r_i - u_i - K\} + f_{i+1}(\min\{Q_i + x_i - u_i, \overline{Q}\})\right\}. \quad (1)$$

We restrict the $u_i$ to lie in the interval $0 \le u_i \le \min\{r_i, Q_i + x_i\}$, which is a natural restriction, since consumption of water $x_i$ together with $Q_i$ at $i^{th}$ period is permanent. Clearly, $f_n(Q_n)$ is easily calculated in view of the fact that best value for $u_n$ is $\overline{u}_n = \min\{r_n, Q_n + x_n, \overline{Q}\}$.

Now consider a cascade system of $\ell$ hydroelectric power stations when the water reservoir is available only for the upper station and between stations $j-1$ and $j$ ($j = 2, \ell$). There is an additional inflow $\xi_i^j$ which is used case of a system of cascaded hydroelectric power stations. Let us suppose that the $j^{th}$ station can pass $I^j$ amount of water. Then this station at the $i$-th period can produce $N_i^j$ energy from the $\min\{I^j, u_i + \xi_i^1 + \xi_i^2 + \cdots + \xi_i^j\}$ amount of water. Thus, we obtain the following recurrent relation:

$$f_i(Q_i) = \min_{\substack{0 \le u_i \le Q_i + x_i \\ \sum_{j=1}^{\ell} N_i^j \le r_i}}\left\{c \cdot \min\left\{K, r_i - \sum_{j=1}^{\ell} N_i^j\right\} + p \cdot \max\left\{0, r_i - \sum_{j=1}^{\ell} N_i^j - K\right\} + f_{i+1}(\min\{Q_i + x_i - u_i, \overline{Q}\})\right\}. \quad (2)$$

As noted above, $f_i(Q_i)$ ($i = 1, \ldots, n-1$) represents the minimum total expenditures in periods from $i$ up to $n$ during the planning horizon in which optimal control of the cascade system of hydroelectric plants is being analyzed when there is $Q_i$ amount of water in a water reservoir at the beginning of $i^{th}$ period. After $\overline{u}_n$ is calculated, $f_n(Q_n)$ can be calculated.

It is possible to write the similar equations for cases when in the cascade there are also other water reservoirs. Generally, when each hydroelectric power station has a water reservoir, we write

$$f_i(Q_i^1, Q_i^2, \ldots, Q_i^\ell) =$$

$$= \min_{\substack{0 \leq u_i^j \leq Q_i^j + \xi_i^j + u_i^{j-1} \\ \sum_{j=1}^{\ell} u_i^j \leq r_i}} \left\{ c \cdot \min\left\{K, r_i - \sum_{j=1}^{\ell} u_i^j\right\} + p \cdot \max\left\{0, r_i - \sum_{j=1}^{\ell} u_i^j - K\right\} + f_{i+1}(Q_{i+1}^1, Q_{i+1}^2, \ldots, Q_{i+1}^\ell) \right\}, \quad i = 1, n-1, \quad (3)$$

where $Q_{i+1}^j = \min\{\overline{Q}^j, Q_i^j + u_i^{j-1} - u_i^j + \xi_i^j\}$. In order to simplify the calculations instead of restriction $\sum_{j=1}^{l} u_i^j \leq r_i$ in braces we can employ a multiplicative the "penalty" term, $\gamma > 0$, so that we have $\gamma \max\left\{0, \sum_{j=1}^{l} u_i^j - r_i\right\}$. When there is possibility to sale the energy, we add the same term with negative coefficient, $\gamma = -a_i$, where $a_i (0 \leq a_i < p)$ is the price of the unit electric energy at the $i^{th}$ period.

Obviously, in fact the relation (3) is suitable at most for the case $\ell = 3$ since for $\ell \geq 4$ the complexity of the problem makes its solution intractable. Similar difficulty arises in the general problem of optimal control of water resources in power system. To avoid this "curse of dimensionality" at least partially, we can introduce some "reasonable" assumptions for the policy of control. For example, we can apply one controlling parameter $u_i$ for $i^{th}$ period (total hydroelectric power produced) and divide its value on stations proportionally with their potential possibilities $Q_i^j + x_i^j$, $j = 1, 2, \ldots, l$.

The recurrent relation for the general case looks like

$$f_i(Q_i^1, Q_i^2, \ldots, Q_i^m) = \min_{0 \leq u_i \leq M} \left\{ c \cdot \min\{K, r_i - u_i\} + p \cdot \max\{0, r_i - u_i - K\} - a \cdot \max\{0, u_i - r_i\} + f_{i+1}(Q_{i+1}^1, Q_{i+1}^2, \ldots, Q_{i+1}^\ell) \right\}, \quad i = 1, n-1, \quad (4)$$

where $M = \sum_{j=1}^{m}(Q_i^j + x_i^j)$, $Q_{i+1}^j = Q_i^j + x_i^j - \frac{u_i}{M} \cdot (Q_i^j + x_i^j)$. At each step it is necessary to perform one-dimensional minimization for $d^1 \cdot d^2 \cdot \ldots \cdot d^m$ sets of initial values of parameters $Q_i^1, Q_i^2, \ldots, Q_i^m$, where $d^j$ specifies the number of possible water levels in $j^{th}$ water reservoir for the standard unit of volume $\Delta q$, $\overline{Q}^j = \Delta q \cdot d^j$.

Note that demand for the hydroelectric power during periods of our analysis is steady and can be predicted. Thus, we shall consider demands as constants. Regarding the flows two hypotheses usually are made – random variables $x_i$ are independent or correlated and they represent Markov chain (simple or complicated). In both cases hypothetical distribution functions are obtained by means of rich statistical data and the relation (1) can be altered in an appropriate way. Respectively $[0; \min\{r_i, Q_i + \overline{\delta}(x_i)\}]$ should be considered as the interval of variation of the parameter $u_i$, where $\overline{\delta}(x_i)$ denotes the upper bound of the interval of fidelity of the random variable $x_i$ with a given significance level.

It must be remarked that for several stations the assumption of independence of variables $x_i$ gives no practical benefit. For one station, in the case when $x_i$'s are independent random variables as in [2, 3], we have:

$$f_i(Q_i) = \min_{u_i} \left\{ c \cdot \min\{K, r_i - u_i\} + p \cdot \max\{0, r_i - u_i - K\} - a \cdot \max\{0, u_i - r_i\} + \int_0^\infty f_{i+1}(\min\{\overline{Q}, Q_i + x_i - u_i\}) \phi_i(x_i) dx_i \right\}, \quad (1')$$

where $\varphi_i(x_i)$ is a density function of water-flow in the $i^{th}$ period.

The case, where $x_i$'s represent a Markov chain as in [1] yields:

$$f_i(Q_i, x_{i-1}) = \min_{u_i} \left\{ c \cdot \min\{K, r_i - u_i\} + p \cdot \max\{0, r_i - u_i - K\} - a \cdot \max\{0, u_i - r_i\} + \int_0^\infty f_{i+1}(\min\{\overline{Q}, Q_i + x_i - u_i\}, x_i) g_i(x_i \mid x_{i-1}) dx_i \right\}, \quad (1'')$$

where $g_i(x_i \mid x_{i-1})$ is a conditional density function. Here on each step it is necessary to perform minimization by $u_i$ for every pair $(Q_i, x_{i-1})$.

Model (4) suggested above was applied to the power system of Georgia. The system included one thermal electric station, three hydro electric stations with water reservoirs and twenty five hydroelectric stations without water reservoirs. The data collected 20 years ago was used. Results are close to the results obtained by J. Giorgobiani, M. Nachkebia, and A. Toronjadze earlier by means of model given in [4]. Thus, the simpler penalization model which avoids the complexity due to higher dimensionality (more hydroelectric power stations and reservoirs in the system to be analyzed) is adequate enough in practice. The simpler model required 60= 5x3x4 one-dimensional minimizations at each step.